\documentclass[11pt,reqno,a4paper]{amsart}
\usepackage{amssymb}
\input amssym.def
\usepackage{amsmath,amsfonts,hyperref,xcolor,mathtools}
\usepackage{amscd,youngtab}
\usepackage[mathscr]{eucal}
\usepackage{enumitem}
\usepackage{orcidlink}
\allowdisplaybreaks
\numberwithin{equation}{section}

\hypersetup{colorlinks=true,citecolor={purple},linkcolor={teal},urlcolor={violet}}

\theoremstyle{plain}
\newtheorem{theorem}{Theorem}[section]
\newtheorem{corollary}[theorem]{Corollary}
\newtheorem{lemma}[theorem]{Lemma}

\theoremstyle{definition}

\newtheorem{remark}[theorem]{Remark}

\makeatletter
\@namedef{subjclassname@2020}{%
\textup{2020} Mathematics Subject Classification}
\makeatother

\subjclass[2020]{11P83, 05A17, 11F27}

\keywords{Congruences, $4$-core partition pairs, $q$-series, dissection formulas}

\begin{document}
\title[Linear identities for partition pairs with $4$-cores]{Linear identities for partition pairs with $4$-cores}
\author[Russelle Guadalupe]{Russelle Guadalupe\orcidlink{0009-0001-8974-4502}}
\address{Institute of Mathematics, University of the Philippines Diliman\\
Quezon City 1101, Philippines}
\email{rguadalupe@math.upd.edu.ph}

\begin{abstract}
We determine an infinite family of linear identities for the number $A_4(n)$ of partition pairs of $n$ with $4$-cores by employing elementary $q$-series techniques and certain $3$-dissection formulas. We then discover an infinite family of congruences for $A_4(n)$ as a consequence of these linear identities.
\end{abstract}

\maketitle

\section{Introduction}\label{sec1}

For complex numbers $a$ and $q$ with $|q| < 1$, define $(a;q)_\infty := \prod_{n=0}^\infty (1-aq^n)$ and $f_m := (q^m;q^m)_\infty$ for integers $m\geq 1$. Recall that a partition $\lambda$ of a positive integer $n$ is a nonincreasing sequence $\lambda_1,\ldots,\lambda_m$ of positive integers, known as the parts of $\lambda$, whose sum is $n$. We describe given partition $\lambda = (\lambda_1,\ldots,\lambda_m)$ of $n$ via its Ferrers--Young diagram consisting of $m$ rows of left-justified boxes, where row $k$ has $\lambda_k$ boxes. We define the hook length of a box in row $i$ and column $j$ as the total number of boxes directly below and to the right of that box, including the box itself. The following figure shows the Ferrers--Young diagram for the partition $(5,3,1,1)$ of $10$ and its corresponding hook lengths.
\begin{align*}
\young(85421,521,2,1)
\end{align*}

We say that a partition of $n$ is $t$-core if none of its hook lengths is divisible by $t$ for some positive integer $t$. For example, the above figure reveals that the partition $(5,3,1,1)$ of $10$ is $t$-core for $t\in \{3,6,7\}$ and $t\geq 9$.

The notion of $t$-core partitions was introduced in relation with the study of $t$-modular representations of the symmetric group; we refer the interested reader to the book \cite{jamker} for a detailed description of $t$-core partitions and its application to these representations. Let $c_t(n)$ be the number of partitions of $n$ that are $t$-cores. Granville and Ono \cite{grano} showed that $c_t(n) > 0$ for any $t\geq 4$ and all positive integers $n$. Garvan, Kim, and Stanton \cite{gks} used combinatorial arguments to prove that the generating function for $c_t(n)$ with $c_t(0) := 1$ is 
\begin{align*}
\sum_{n=0}^\infty c_t(n)q^n = \dfrac{f_t^t}{f_1}.
\end{align*}

Several authors explored arithmetic properties of $c_t(n)$ via classical $q$-series techniques and modular forms. Garvan \cite{garv} found Ramanujan-type congruences for $c_p(n)$ for some suitable primes $p$. Hirschhorn and Sellers \cite{hirscsela,hirscselb} proved congruences modulo $2$ and $4$ and linear identities for $c_t(n)$ such as 
\begin{align*}
c_4\left(3^{2k+1}n+\dfrac{5\cdot 3^{2k}-5}{8}\right) &= 3^kc_4(3n),\\
c_4\left(3^{2k+1}n+\dfrac{13\cdot 3^{2k}-5}{8}\right) &= (2\cdot 3^k-1)c_4(3n+1),\\
c_4\left(3^{2k+2}n+\dfrac{7\cdot 3^{2k+1}-5}{8}\right) &= \dfrac{3^{k+1}-1}{2}c_4(9n+2),\\
c_4\left(3^{2k+2}n+\dfrac{23\cdot 3^{2k+1}-5}{8}\right) &= \dfrac{3^{k+1}-1}{2}c_4(9n+8)
\end{align*}
for all $n\geq 0$ and $k\geq 1$. Ono and Sze \cite{onosze} obtained a formula for $c_4(n)$ in terms of the class number of a certain imaginary quadratic field and settled other linear identities conjectured by Hirschhorn and Sellers \cite{hirscsela}. Baruah and Berndt \cite{barber} found linear identities for $c_3(n)$ and $c_5(n)$ by employing Ramanujan's classical modular equations of third and fifth degrees, respectively. Kim \cite{kim} demonstrated a unified way of discovering linear identities for $c_p(n)$ for primes $p\geq 5$ by applying Hecke operators. Baruah and Nath \cite{barnata} deduced infinite families of linear identities for $a_4(n)$ analogous to the above identities of Hirschhorn and Sellers \cite{hirscsela}. We refer the reader to the survey paper \cite{chokns} for more results on the arithmetic properties of $t$-core partitions.

We define a partition pair of $n$ with $t$-cores by a pair of partitions $(\lambda, \mu)$ such that the sum of all parts of $\lambda$ and $\mu$ is $n$ and both $\lambda$ and $\mu$ are $t$-cores. We see that the generating function for the number $A_t(n)$ of partition pairs of $n$ with $t$-cores is given by
\begin{align*}
\sum_{n=0}^\infty A_t(n)q^n := \dfrac{f_t^{2t}}{f_1^2}.
\end{align*} 

Similar to that of $c_t(n)$, various authors studied arithmetic properties of $A_t(n)$. Lin \cite{lin} found congruences modulo $4, 5, 7$, and $8$ for $A_3(n)$. Baruah and Nath \cite{barnatb} derived infinite families of linear identities for $A_3(n)$ via Ramanujan's theta function identities, generalizing one of the results of Lin \cite{lin}. Yao \cite{yao} proved congruences modulo $3$ and $9$ for $A_3(n)$, and Xia \cite{xia} did the same for congruences modulo $(4^m-1)/3$ for $A_3(n)$. Dasappa \cite{dasa} extended the results of Saikia and Boruah \cite{saiborb} by proving infinite family of congruences modulo powers of $5$ for $A_5(n)$. Saikia and Boruah \cite{saibora} showed congruences modulo $2$ for $A_4(n)$. Recently, the author \cite{guad} determined an infinite family of linear identities for $A_5(n)$ and deduced new congruences for $A_5(n)$.

We explore in this paper arithmetic properties of $A_4(n)$ by utilizing elementary $q$-series manipulations. Our main result reveals the following generating function for integers $k\geq 0$.

\begin{theorem}\label{thm11}
For integers $k\geq 0$, we have
\begin{align}
\sum_{n=0}^\infty A_4\left(3^{2k+1}n+\dfrac{3^{2k+2}-5}{4}\right)q^n &= \dfrac{3^{2k+2}-1}{4}\cdot\dfrac{f_2^{12}}{f_1^6}\nonumber\\
&+\dfrac{27(3^{2k}-1)(3^{2k+2}-1)}{320}\left(3\dfrac{f_2^2f_3^{10}}{f_1^4f_6^2}+4q\dfrac{f_3f_6^7}{f_1f_2}\right)+3^{2k+2}\cdot\dfrac{q^3f_{12}^8}{f_3^2}.\label{eq11}
\end{align}
\end{theorem}

A quick comparison shows that the above generating function is more complicated than that of \cite[Theorem 1.1]{guad}. A consequence of Theorem \ref{thm11} is the following infinite family of linear identities for $A_4(n)$ and an infinite family of congruences for $A_4(n)$.

\begin{theorem}\label{thm12}
For all integers $n\geq 0$ with $3\nmid n$ and $k\geq 0$, we have 
\begin{align}\label{eq12}
A_4\left(3^{2k+1}n+\dfrac{3^{2k+2}-5}{4}\right)= \dfrac{3^{4k+4}-1}{80}A_4(3n+1).
\end{align}
Consequently, for all integers $n\geq 0$ with $3\nmid n$ and $k\geq 0$, we have 
\begin{align}\label{eq13}
A_4\left(3^{2k+1}n+\dfrac{3^{2k+2}-5}{4}\right) \equiv 0\pmod{\dfrac{3(3^{4k+4}-1)}{40}}.
\end{align}
\end{theorem}

We organize the rest of the paper as follows. We present in Section \ref{sec2} some theta function identities needed to prove Theorem \ref{thm11}, including $3$-dissection formulas and identities involving the Ramanujan's cubic continued fraction $C(q)$ and the continued fraction $h(q)$ of level $12$. In Section \ref{sec3}, we apply these identities to find the generating functions for $A_4(3n+1)$ and $A_4(27n+19)$ in Section \ref{sec3}. With the help of these generating functions, we finally demonstrate in Section \ref{sec4} Theorems \ref{thm11} and \ref{thm12}. We also prove another linear identity for $A_4(n)$ similar to (\ref{eq12}) and give some comments. We have performed our calculations via \textit{Mathematica}. 

\section{Some necessary identities}\label{sec2}

We provide in this section important identities that will assist us in showing Theorems \ref{thm11} and \ref{thm12}. These include certain $3$-dissections formulas and some identities involving the Ramanujan's cubic continued fraction $C(q)$ given by  
\begin{align*}
	C(q) := \dfrac{1}{1+\dfrac{q+q^2}{1+\dfrac{q^2+q^4}{1+\dfrac{q^3+q^6}{1+\cdots}}}} = \dfrac{f_1f_6^3}{f_2f_3^3}
\end{align*}
and the continued fraction $h(q)$ of order $12$ defined by Mahadeva Naika, B. N. Dharmendra, and Shivashankara \cite{mnds}, and is given by 
\begin{align*}
	h(q) &:= \dfrac{q(1-q)}{1-q^3+\dfrac{q^3(1-q^2)(1-q^4)}{(1-q^3)(1+q^6)+\dfrac{q^3(1-q^8)(1-q^{10})}{(1-q^3)(1+q^{12})+\cdots}}}\\
	&= q\prod_{n=1}^\infty \dfrac{(1-q^{12n-1})(1-q^{12n-11})}{(1-q^{12n-5})(1-q^{12n-7})}.
\end{align*}
We begin with the following $3$-dissection formulas \cite[p. 195, (22.6.13)]{hirsc}, \cite[p. 68]{hirscselc}
\begin{align}
	\dfrac{f_2^2}{f_1} &= \dfrac{f_6f_9^2}{f_3f_{18}}+q\dfrac{f_{18}^2}{f_9},\label{eq21}\\
	\dfrac{f_2}{f_1^2} &= \dfrac{f_6^4f_9^6}{f_3^8f_{18}^3}+2q\dfrac{f_6^3f_9^3}{f_3^7}+4q^2\dfrac{f_6^2f_{18}^3}{f_3^6}.\label{eq22}
\end{align}

We next require the following identities involving $C(q)$ \cite[Theorem 6.9, Theorem 6.11]{coop}, \cite[Theorem 1, Theorem 2]{chan}:
\begin{align}
	\dfrac{f_2^8f_3^4}{qf_1^4f_6^8} &= \dfrac{1}{qC^3(q)}+1,\label{eq23}\\
	\dfrac{f_1^4f_2^4}{qf_3^4f_6^4} &= \dfrac{1}{qC^3(q)}-7-8qC^3(q),\label{eq24}\\
	f_1f_2 &= f_9f_{18}\left(\dfrac{1}{C(q^3)}-q-2q^2C(q^3)\right),\label{eq25}\\
	\dfrac{1}{f_1f_2} &= \dfrac{f_9^3f_{18}^3}{f_3^4f_6^4}\left(\dfrac{1}{C^2(q^3)}+\dfrac{q}{C(q^3)}+3q^2-2q^3C(q^3)+4q^4C^2(q^3)\right). \label{eq26}
\end{align}
We remark that (\ref{eq25}) follows from replacing $q$ with $q^3$ in \cite[Theorem 2]{chan} and (\ref{eq26}) appears in the proof of \cite[Theorem 1]{chan} (more precisely, it appears as equation (13) in \cite[p. 677]{chan}).

We also need the following identities involving $h(q)$ \cite[Theorem 12.2]{coop}:
\begin{align}
	\dfrac{f_3^3f_4}{qf_1f_{12}^3} &= \dfrac{1}{h(q)}+h(q),\label{eq27}\\
	\dfrac{f_4^4f_6^2}{qf_2^2f_{12}^4} &= \dfrac{1}{h(q)}-1+h(q),\label{eq28}\\
	\dfrac{f_1f_4^2f_6^9}{qf_2^3f_3^3f_{12}^6} &= \dfrac{1}{h(q)}-2+h(q),\label{eq29}\\
	\dfrac{f_2^7f_3}{qf_1^3f_4^2f_6f_{12}^2} &= \dfrac{1}{h(q)}+2+h(q).\label{eq210}
\end{align}

We now use (\ref{eq23})--(\ref{eq210}) to deduce the following auxiliary theta function identities.

\begin{lemma}\label{lem21}
We have the identities
\begin{align}
\dfrac{f_2^7f_3^5}{f_1^5f_6} &= \dfrac{f_2^2f_3^{10}}{f_1^4f_6^2}+q\dfrac{f_3f_6^7}{f_1f_2}, \label{eq211}\\
\dfrac{f_1^3f_2^3f_6^3}{f_3^3} &=-8\dfrac{f_2^{12}}{f_1^6}+9\left(\dfrac{f_2^2f_3^{10}}{f_1^4f_6^2}+q\dfrac{f_3f_6^7}{f_1f_2}\right).\label{eq212}
\end{align}
\end{lemma}

\begin{proof}
Multiplying both sides of (\ref{eq23}) by $qf_3f_6^7/(f_1f_2)$ yields (\ref{eq211}). On the other hand, dividing both sides of (\ref{eq212}) by $qf_3f_6^7/(f_1f_2)$ which is the last term of this
equation, it suffices to show that
\begin{align}\label{eq213}
\dfrac{f_1^4f_2^4}{qf_3^4f_6^4}+8\dfrac{f_2^{13}}{qf_1^5f_3f_6^7} = 9\dfrac{f_2^3f_3^9}{qf_1^3f_6^9}+9.
\end{align}
Observe that from (\ref{eq23}), we have
\begin{align}\label{eq214}
\dfrac{f_2^{13}}{qf_1^5f_3f_6^7} = \left(\dfrac{f_2^8f_3^4}{qf_1^4f_6^8}\right)^2\cdot \dfrac{qf_1^3f_6^9}{f_2^3f_3^9} = \left(\dfrac{1}{qC^3(q)}+1\right)^2\cdot qC^3(q).
\end{align}
Combining (\ref{eq24}) and (\ref{eq214}) leads us to
\begin{align*}
\dfrac{f_1^4f_2^4}{qf_3^4f_6^4}+8\dfrac{f_2^{13}}{qf_1^5f_3f_6^7} &= \dfrac{1}{qC^3(q)}-7-8qC^3(q)+8\left(\dfrac{1}{qC^3(q)}+2+qC^3(q)\right)\\
&= \dfrac{9}{qC^3(q)}+9=9\dfrac{f_2^3f_3^9}{qf_1^3f_6^9}+9,
\end{align*}
which is precisely (\ref{eq213}). This completes the proof of (\ref{eq212}).
\end{proof}

\begin{lemma}\label{lem22}
We have the identity
\begin{align*}
\dfrac{f_2^{15}f_3^2f_{12}}{qf_1^6f_4^3f_6^9}+4q^2\dfrac{f_2^3f_{12}^9}{f_4^3f_6^9}=\dfrac{f_3^3f_4}{qf_1f_{12}^3}+2+3\dfrac{f_2^3f_3^6f_{12}^3}{f_1^2f_4f_6^9}+4q\dfrac{f_2^3f_3^3f_{12}^6}{f_1f_4^2f_6^9}.
\end{align*}
\end{lemma}

\begin{proof}
We recast the given identity as 
\begin{align}\label{eq215}
\dfrac{f_2^{15}f_3^2f_{12}}{qf_1^6f_4^3f_6^9}-3\dfrac{f_2^3f_3^6f_{12}^3}{f_1^2f_4f_6^9}-4q\dfrac{f_2^3f_3^3f_{12}^6}{f_1f_4^2f_6^9}+4q^2\dfrac{f_2^3f_{12}^9}{f_4^3f_6^9}=\dfrac{f_3^3f_4}{qf_1f_{12}^3}+2.
\end{align}
Letting $h := h(q)$, we write each summand on the left-hand side of (\ref{eq215}) using (\ref{eq27})--(\ref{eq210}). We find that
\begin{align}
\dfrac{f_2^{15}f_3^2f_{12}}{qf_1^6f_4^3f_6^9} &= \dfrac{f_4^4f_6^2}{qf_2^2f_{12}^4}\cdot\dfrac{qf_1f_{12}^3}{f_3^3f_4}\left(\dfrac{f_2^7f_3}{qf_1^3f_4^2f_6f_{12}^2}\right)^2\dfrac{qf_2^3f_3^3f_{12}^6}{f_1f_4^2f_6^9}\nonumber\\
&=\left(\dfrac{1}{h}-1+h\right)\left(\dfrac{1}{h}+h\right)^{-1}\left(\dfrac{1}{h}+2+h\right)^2\left(\dfrac{1}{h}-2+h\right)^{-1}\nonumber\\
&=\dfrac{(1+h)^4(1-h+h^2)}{h(1-h)^2(1+h^2)}\label{eq216},\\
\dfrac{f_2^3f_3^6f_{12}^3}{f_1^2f_4f_6^9} &= \dfrac{f_3^3f_4}{qf_1f_{12}^3}\cdot \dfrac{qf_2^3f_3^3f_{12}^6}{f_1f_4^2f_6^9}= \left(\dfrac{1}{h}+h\right)\left(\dfrac{1}{h}-2+h\right)^{-1}\nonumber\\
&= \dfrac{1+h^2}{(1-h)^2}\label{eq217s},\\
q\dfrac{f_2^3f_3^3f_{12}^6}{f_1f_4^2f_6^9} &= \left(\dfrac{1}{h}-2+h\right)^{-1} = \dfrac{h}{(1-h)^2}\label{eq218},\\
q^2\dfrac{f_2^3f_{12}^9}{f_4^3f_6^9} &= \dfrac{qf_1f_{12}^3}{f_3^3f_4}\cdot \dfrac{qf_2^3f_3^3f_{12}^6}{f_1f_4^2f_6^9}=\left(\dfrac{1}{h}+h\right)^{-1}\left(\dfrac{1}{h}-2+h\right)^{-1} \nonumber\\
&=\dfrac{h^2}{(1-h)^2(1+h^2)}. \label{eq219}
\end{align}

We infer from (\ref{eq27}) and (\ref{eq216})--(\ref{eq219}) that
\begin{align*}
\dfrac{f_2^{15}f_3^2f_{12}}{qf_1^6f_4^3f_6^9}&-3\dfrac{f_2^3f_3^6f_{12}^3}{f_1^2f_4f_6^9}-4q\dfrac{f_2^3f_3^3f_{12}^6}{f_1f_4^2f_6^9}+4q^2\dfrac{f_2^3f_{12}^9}{f_4^3f_6^9}\\
&= \left(\dfrac{(1+h)^4}{h(1-h)^2}-\dfrac{(1+h)^4}{(1-h)^2(1+h^2)}\right)-\dfrac{3(1+h^2)}{(1-h)^2}-\dfrac{4h}{(1-h)^2}+\dfrac{4h^2}{(1-h)^2(1+h^2)}\\
&= \dfrac{(1+h)^4}{h(1-h)^2}-\dfrac{3h^2+4h+3}{(1-h)^2}+\dfrac{(2h+(1+h)^2)(2h-(1+h)^2)}{(1-h)^2(1+h^2)}\\
&= \dfrac{(1+h)^4}{h(1-h)^2}-\dfrac{3h^2+4h+3}{(1-h)^2}-\dfrac{h^2+4h+1}{(1-h)^2}\\
&= \dfrac{(1+h)^4}{h(1-h)^2}-\dfrac{4(1+h)^2}{(1-h)^2}= \left(\dfrac{1+h}{1-h}\right)^2\left(\dfrac{(1+h)^2}{h}-4\right)\\
&= \left(\dfrac{1+h}{1-h}\right)^2\left(\dfrac{(1-h)^2}{h}\right) = \dfrac{1}{h}+h+2 = \dfrac{f_3^3f_4}{qf_1f_{12}^3}+2,
\end{align*}
which is the right-hand side of (\ref{eq215}). This completes the proof.
\end{proof}

\section{Generating functions for $A_4(3n+1)$ and $A_4(27n+19)$}\label{sec3}

We obtain in this section generating functions for $A_4(3n+1)$ and $A_4(27n+19)$ using the theta function identities proved in Section \ref{sec2}. These will be instrumental to the proofs of Theorem \ref{thm11} and \ref{thm12}, which will be explained in Section \ref{sec4}. We begin with presenting the generating function for $A_4(3n+1)$ and a congruence modulo $6$ for $A_4(n)$ as a byproduct. 

\begin{theorem}\label{thm31}
We have that
\begin{align}\label{eq31}
\sum_{n=0}^\infty A_4(3n+1)q^n = 2\dfrac{f_2^{12}}{f_1^6}+9q^3\dfrac{f_{12}^8}{f_3^2}.
\end{align}
\end{theorem}

\begin{proof}
We apply (\ref{eq21}) twice in the generating function for $A_4(n)$ found in the introduction: once with $q$ replaced with $q^2$, and once in its original form. We obtain
\begin{align}
\sum_{n=0}^\infty A_4(n)q^n &= \dfrac{f_4^8}{f_1^2} = \left(\dfrac{f_4^2}{f_2}\right)^4\left(\dfrac{f_2^2}{f_1}\right)^2\nonumber\\
&= \left(\dfrac{f_{12}f_{18}^2}{f_6f_{36}}+q^2\dfrac{f_{36}^2}{f_{18}}\right)^4\left(\dfrac{f_6f_9^2}{f_3f_{18}}+q\dfrac{f_{18}^2}{f_9}\right)^2.\label{eq32}
\end{align}
We consider the terms involving $q^{3n+1}$ in the expansion of (\ref{eq32}). We see that
\begin{align*}
\sum_{n=0}^\infty A_4(3n+1)q^n &= 2\dfrac{f_3f_4^4f_6^9}{f_1f_2^3f_{12}^4}+4q\dfrac{f_4^3f_6^9}{f_2^3f_3^2f_{12}}+6q\dfrac{f_3^4f_4^2f_{12}^2}{f_1^2}+8q^2\dfrac{f_3f_4f_{12}^5}{f_1}+q^3\dfrac{f_{12}^8}{f_3^2}\\
&=2q\dfrac{f_4^3f_6^9}{f_2^3f_3^2f_{12}}\left(\dfrac{f_3^3f_4}{qf_1f_{12}^3}+2+3\dfrac{f_2^3f_3^6f_{12}^3}{f_1^2f_4f_6^9}+4q\dfrac{f_2^3f_3^3f_{12}^6}{f_1f_4^2f_6^9}\right)+q^3\dfrac{f_{12}^8}{f_3^2}.
\end{align*}
Applying Lemma \ref{lem22} on the expression inside the parentheses in the above equality leads to
\begin{align*}
\sum_{n=0}^\infty A_4(3n+1)q^n &=2q\dfrac{f_4^3f_6^9}{f_2^3f_3^2f_{12}}\left(\dfrac{f_2^{15}f_3^2f_{12}}{qf_1^6f_4^3f_6^9}+4q^2\dfrac{f_2^3f_{12}^9}{f_4^3f_6^9}\right)+q^3\dfrac{f_{12}^8}{f_3^2}\\
&=2\dfrac{f_2^{12}}{f_1^6}+8q^3\dfrac{f_{12}^8}{f_3^2}+q^3\dfrac{f_{12}^8}{f_3^2}\\
&=2\dfrac{f_2^{12}}{f_1^6}+9q^3\dfrac{f_{12}^8}{f_3^2}
\end{align*}
as desired.
\end{proof}

\begin{theorem}\label{thm32}
For all integers $n\geq 0$ with $3\nmid n$, we have $A_4(3n+1) \equiv 0\pmod{6}$.
\end{theorem}

\begin{proof}
Using (\ref{eq21}), we express (\ref{eq31}) as
\begin{align}\label{eq33}
\sum_{n=0}^\infty A_4(3n+1)q^n = 2\left(\dfrac{f_2^2}{f_1}\right)^6+9q^3\dfrac{f_{12}^8}{f_3^2}=2\left(\dfrac{f_6f_9^2}{f_3f_{18}}+q\dfrac{f_{18}^2}{f_9}\right)^6+9q^3\dfrac{f_{12}^8}{f_3^2}.
\end{align}
Extracting the terms involving $q^{3n+1}$ in the expansion of (\ref{eq33}), dividing by $q$, and then replacing $q^3$ with $q$, we get
\begin{align*}
\sum_{n=0}^\infty A_4(9n+4)q^n &= 6\left(2\dfrac{f_2^5f_3^9}{f_1^5f_6^3}+5q\dfrac{f_2^2f_6^6}{f_1^2}\right).
\end{align*}
Similarly, we extract the terms involving $q^{3n+2}$ in the expansion of (\ref{eq33}), divide by $q$, and then replace $q^3$ with $q$, so that
\begin{align*}
\sum_{n=0}^\infty A_4(9n+7)q^n &= 6\left(5\dfrac{f_2^4f_3^6}{f_1^4}+2q\dfrac{f_2f_6^9}{f_1f_3^3}\right).
\end{align*}
Hence, we obtain $A_4(9n+4)\equiv 0\pmod{6}$ and $A_4(9n+7)\equiv 0\pmod{6}$ for all $n\geq 0$ as desired.
\end{proof}

Given a $q$-series $f(q) := \sum_{n = n_0}^\infty a(n)q^n\in \mathbb{Z}[[q]]$, we define the linear maps in $\mathbb{Z}[[q]]$ by
\begin{align*}
U(f(q)) &:= \sum_{n=\lceil n_0/3\rceil}^\infty a(3n)q^n,\\
V(f(q)) &:= \sum_{n=\lceil (n_0-1)/3\rceil}^\infty a(3n+1)q^n.
\end{align*}

The next two lemmas provide the closed-form formulas for the images of certain $q$-series under $U$ and $V$. This will play a key role in determining the generating function for $A_4(27n+19)$ and in showing Theorem \ref{thm11} in Section \ref{sec4}.

\begin{lemma}\label{lem33}
We have the identities 
\begin{align}
U\left(\dfrac{f_2^{12}}{q^3f_1^6}\right) &= \dfrac{f_2^6f_3^{12}}{qf_1^6f_6^6}+20\dfrac{f_2^3f_3^3f_6^3}{f_1^3}+q\dfrac{f_6^{12}}{f_3^6},\label{eq34}\\
U\left(\dfrac{f_2^2f_3^{10}}{q^3f_1^4f_6^2}\right) &= \dfrac{f_2^6f_3^{12}}{qf_1^6f_6^6}+16\dfrac{f_2^3f_3^3f_6^3}{f_1^3},\label{eq35}\\
U\left(\dfrac{f_3f_6^7}{q^2f_1f_2}\right) &=  3\dfrac{f_2^3f_3^3f_6^3}{f_1^3}.\label{eq36}
\end{align}
\end{lemma}

\begin{proof}
We use (\ref{eq21}) to write 
\begin{align}\label{eq37}
\dfrac{f_2^{12}}{q^3f_1^6} = \dfrac{1}{q^3}\left(\dfrac{f_2^2}{f_1}\right)^6=\dfrac{1}{q^3}\left(\dfrac{f_6f_9^2}{f_3f_{18}}+q\dfrac{f_{18}^2}{f_9}\right)^6.
\end{align}
Examining the terms involving $q^{3n}$ in the expansion of (\ref{eq37}), we get (\ref{eq34}).
We again use (\ref{eq21}) to write 
\begin{align}\label{eq38}
\dfrac{f_2^2f_3^{10}}{q^3f_1^4f_6^2} = \dfrac{f_3^{10}}{q^3f_6^2}\left(\dfrac{f_2}{f_1^2}\right)^2=\dfrac{f_3^{10}}{q^3f_6^2}\left(\dfrac{f_6^4f_9^6}{f_3^8f_{18}^3}+2q\dfrac{f_6^3f_9^3}{f_3^7}+4q^2\dfrac{f_6^2f_{18}^3}{f_3^6}\right)^2.
\end{align}
Getting the terms involving $q^{3n}$ in the expansion of (\ref{eq38}), we obtain (\ref{eq35}).
We finally use (\ref{eq26}) so that
\begin{align}\label{eq39}
\dfrac{f_3f_6^7}{q^2f_1f_2}=\dfrac{f_3f_6^7}{q^2}\cdot \dfrac{f_9^3f_{18}^3}{f_3^4f_6^4}\left(\dfrac{1}{C^2(q^3)}+\dfrac{q}{C(q^3)}+3q^2-2q^3C(q^3)+4q^4C^2(q^3)\right).
\end{align}
Examining the terms involving $q^{3n}$ in the expansion of (\ref{eq39}), we arrive at (\ref{eq36}).
\end{proof}

\begin{lemma}\label{lem34}
We have the identities 
\begin{align}
V\left(\dfrac{f_2^6f_3^{12}}{qf_1^6f_6^6}\right) &= 21\dfrac{f_2^2f_3^{10}}{f_1^4f_6^2}+48q\dfrac{f_3f_6^7}{f_1f_2},\label{eq310}\\
V\left(\dfrac{f_2^3f_3^3f_6^3}{f_1^3}\right) &= 3\dfrac{f_2^2f_3^{10}}{f_1^4f_6^2}+3q\dfrac{f_3f_6^7}{f_1f_2}.\label{eq311}
\end{align}
\end{lemma}

\begin{proof}
We know from (\ref{eq21}) and (\ref{eq22}) that
\begin{align}\label{eq312}
\dfrac{f_2^6f_3^{12}}{qf_1^6f_6^6} &= \dfrac{f_3^{12}}{qf_6^6}\left(\dfrac{f_2^2}{f_1}\right)^2\left(\dfrac{f_2}{f_1^2}\right)^2\nonumber\\
&=\dfrac{f_3^{12}}{qf_6^6}\left(\dfrac{f_6f_9^2}{f_3f_{18}}+q\dfrac{f_{18}^2}{f_9}\right)^2\left(\dfrac{f_6^4f_9^6}{f_3^8f_{18}^3}+2q\dfrac{f_6^3f_9^3}{f_3^7}+4q^2\dfrac{f_6^2f_{18}^3}{f_3^6}\right)^2.
\end{align}
Considering the terms involving $q^{3n+1}$ in the expansion of (\ref{eq312}) yields (\ref{eq310}). We also know that
\begin{align}\label{eq313}
\dfrac{f_2^3f_3^3f_6^3}{f_1^3} &= f_3^3f_6^3\left(\dfrac{f_2^2}{f_1}\right)\left(\dfrac{f_2}{f_1^2}\right)\nonumber\\
&=f_3^3f_6^3\left(\dfrac{f_6f_9^2}{f_3f_{18}}+q\dfrac{f_{18}^2}{f_9}\right)\left(\dfrac{f_6^4f_9^6}{f_3^8f_{18}^3}+2q\dfrac{f_6^3f_9^3}{f_3^7}+4q^2\dfrac{f_6^2f_{18}^3}{f_3^6}\right).
\end{align}
Considering the terms involving $q^{3n+1}$ in the expansion of (\ref{eq313}), we see that
\begin{align*}
V\left(\dfrac{f_2^3f_3^3f_6^3}{f_1^3}\right) = 3\dfrac{f_2^7f_3^5}{f_1^5f_6}.
\end{align*}
Thus, (\ref{eq311}) follows from the above expression and (\ref{eq211}).
\end{proof}

\begin{theorem}\label{thm35}
We have that
\begin{align}\label{eq314}
\sum_{n=0}^\infty A_4(27n+19)q^n = 20\dfrac{f_2^{12}}{f_1^6}+54\left(3\dfrac{f_2^2f_3^{10}}{f_1^4f_6^2}+4q\dfrac{f_3f_6^7}{f_1f_2}\right)+81q^3\dfrac{f_{12}^8}{f_3^2}.
\end{align}
\end{theorem}

\begin{proof}
Dividing both sides of (\ref{eq31}) by $q^3$ yields
\begin{align}\label{eq315}
\sum_{n=-3}^\infty A_4(3n+10)q^n = 2\dfrac{f_2^{12}}{q^3f_1^6}+9\dfrac{f_{12}^8}{f_3^2}.
\end{align}
Applying $U$ on both sides of (\ref{eq315}) and employing Lemma \ref{lem33}, we find that
\begin{align}\label{eq316}
\sum_{n=-1}^\infty A_4(9n+10)q^n &= 2U\left(\dfrac{f_2^{12}}{q^3f_1^6}\right)+9\dfrac{f_4^8}{f_1^2}\nonumber\\
&= 2\dfrac{f_2^6f_3^{12}}{qf_1^6f_6^6}+40\dfrac{f_2^3f_3^3f_6^3}{f_1^3}+2q\dfrac{f_6^{12}}{f_3^6}+9\dfrac{f_4^8}{f_1^2}.
\end{align} 
Recall that $f_4^8/f_1^2$ is the generating function for $A_4(n)$. We now apply $V$ on both sides of (\ref{eq316}). In view of Theorem \ref{thm31}, Lemma \ref{lem34}, and the fact that $V\left(\sum_{n=0}^\infty A_4(n)q^n\right) = \sum_{n=0}^\infty A_4(3n+1)q^n$, we deduce that
\begin{align*}
\sum_{n=0}^\infty A_4(27n+19)q^n & = 2V\left(\dfrac{f_2^6f_3^{12}}{qf_1^6f_6^6}\right)+40V\left(\dfrac{f_2^3f_3^3f_6^3}{f_1^3}\right)+2\dfrac{f_2^{12}}{f_1^6}+9\sum_{n=0}^\infty A_4(3n+1)q^n\\
&= 2\left(21\dfrac{f_2^2f_3^{10}}{f_1^4f_6^2}+48q\dfrac{f_3f_6^7}{f_1f_2}\right)+40\left(3\dfrac{f_2^2f_3^{10}}{f_1^4f_6^2}+3q\dfrac{f_3f_6^7}{f_1f_2}\right)\\
&+ 2\dfrac{f_2^{12}}{f_1^6}+9\left(2\dfrac{f_2^{12}}{f_1^6}+9q^3\dfrac{f_{12}^8}{f_3^2}\right),\\
&= 20\dfrac{f_2^{12}}{f_1^6}+162\dfrac{f_2^2f_3^{10}}{f_1^4f_6^2}+216q\dfrac{f_3f_6^7}{f_1f_2}+81q^3\dfrac{f_{12}^8}{f_3^2},
\end{align*}
which is equivalent to (\ref{eq314}).
\end{proof}

We now employ Theorem \ref{thm35} to derive a simple relation between $A_4(3n+1)$ and $A_4(27n+19)$.

\begin{theorem}\label{thm36}
For all integers $n\geq 0$ with $3\nmid n$, we have $A_4(27n+19) = 82A_4(3n+1)$.
\end{theorem}

\begin{proof}
In view of (\ref{eq212}) and Theorems \ref{thm31} and \ref{thm35}, we consider the expansion of 
\begin{align}
\sum_{n=0}^\infty &(A_4(27n+19)-82A_4(3n+1))q^n \nonumber\\
&= 20\dfrac{f_2^{12}}{f_1^6}+54\left(3\dfrac{f_2^2f_3^{10}}{f_1^4f_6^2}+4q\dfrac{f_3f_6^7}{f_1f_2}\right)+81q^3\dfrac{f_{12}^8}{f_3^2}-82\left(2\dfrac{f_2^{12}}{f_1^6}+9q^3\dfrac{f_{12}^8}{f_3^2}\right)\nonumber\\
&= 18\left(-8\dfrac{f_2^{12}}{f_1^6}+9\dfrac{f_2^2f_3^{10}}{f_1^4f_6^2}+12q\dfrac{f_3f_6^7}{f_1f_2}\right) - 657q^3\dfrac{f_{12}^8}{f_3^2}\nonumber\\
&= 18\left(\dfrac{f_1^3f_2^3f_6^3}{f_3^3}-9q\dfrac{f_3f_6^7}{f_1f_2}+12q\dfrac{f_3f_6^7}{f_1f_2}\right) - 657q^3\dfrac{f_{12}^8}{f_3^2}\nonumber\\
&= 18\left(\dfrac{f_1^3f_2^3f_6^3}{f_3^3}+3q\dfrac{f_3f_6^7}{f_1f_2}\right) - 657q^3\dfrac{f_{12}^8}{f_3^2}.\label{eq317}
\end{align}
Utilizing (\ref{eq25}) and (\ref{eq26}), we see that
\begin{align}
\dfrac{f_1^3f_2^3f_6^3}{f_3^3}&+3q\dfrac{f_3f_6^7}{f_1f_2}\nonumber\\
&=\dfrac{f_6^3}{f_3^3}\cdot f_9^3f_{18}^3\left(\dfrac{1}{C(q^3)}-q-2q^2C(q^3)\right)^3+3qf_3f_6^7\cdot \dfrac{f_9^3f_{18}^3}{f_3^4f_6^4}\left(\dfrac{1}{C^2(q^3)}+\dfrac{q}{C(q^3)}+3q^2\right.\nonumber\\
&\left.-2q^3C(q^3)+4q^4C^2(q^3)\right)\nonumber\\
&= \dfrac{f_6^3f_9^3f_{18}^3}{f_3^3}\left(\dfrac{1}{C^3(q^3)}+20q^3-8q^6C^3(q^3)\right).\label{eq318}
\end{align}
We infer from (\ref{eq317}) and (\ref{eq318}) that the right-hand side of (\ref{eq317}) is a $q$-series whose terms are of the form $q^{3n}$ only. Thus, the coefficients of the terms involving $q^{3n+1}$ and $q^{3n+2}$ in (\ref{eq317}) are all zero, arriving at the desired identity.
\end{proof}

\begin{corollary}\label{cor37}
For all integers $n\geq 0$ with $3\nmid n$, we have $A_4(27n+19)\equiv 0\pmod{492}$.
\end{corollary}

\begin{proof}
This immediately follows from Theorems \ref{thm32} and \ref{thm36}.
\end{proof}

\section{Proofs of Theorems \ref{thm11} and \ref{thm12}}\label{sec4}

We demonstrate in this section Theorems \ref{thm11} and \ref{thm12} by employing the generating functions for $A_4(3n+1)$ and $A_4(27n+19)$ and the $q$-series identities involving the linear maps $U$ and $V$ determined in Section \ref{sec3}. We remark that Theorem \ref{thm32} and Corollary \ref{cor37} are particular cases of (\ref{eq13}) for $k=0$ and $k=1$, respectively.

\begin{proof}[Proof of Theorem \ref{thm11}]
We proceed by induction on $k$. Appealing to Theorems \ref{thm31} and \ref{thm35}, we see that (\ref{eq11}) holds for $k\in \{0,1\}$. Suppose now that (\ref{eq11}) holds for some $k\geq 2$. Dividing both sides of (\ref{eq11}) by $q^3$ yields
\begin{align}
	\sum_{n=-3}^\infty &A_4\left(3^{2k+1}n+\dfrac{5(3^{2k+2}-1)}{4}\right)q^n \nonumber\\ &=\dfrac{3^{2k+2}-1}{4}\cdot\dfrac{f_2^{12}}{q^3f_1^6}+\dfrac{27(3^{2k}-1)(3^{2k+2}-1)}{320}\left(3\dfrac{f_2^2f_3^{10}}{q^3f_1^4f_6^2}+4\dfrac{f_3f_6^7}{q^2f_1f_2}\right)+3^{2k+2}\cdot\dfrac{f_{12}^8}{f_3^2}.\label{eq41}
\end{align}
We apply $U$ on both sides of (\ref{eq41}) and use Lemma \ref{lem33}, obtaining
\begin{align}\label{eq42}
	\sum_{n=-1}^\infty &A_4\left(3^{2k+2}n+\dfrac{5(3^{2k+2}-1)}{4}\right)q^n \nonumber\\
	&= \dfrac{3^{2k+2}-1}{4}U\left(\dfrac{f_2^{12}}{q^3f_1^6}\right)+\dfrac{27(3^{2k}-1)(3^{2k+2}-1)}{320}\left[3U\left(\dfrac{f_2^2f_3^{10}}{q^3f_1^4f_6^2}\right)+4U\left(\dfrac{f_3f_6^7}{q^2f_1f_2}\right)\right]\nonumber\\
	&+3^{2k+2}\cdot\dfrac{f_4^8}{f_1^2}\nonumber\\
	&= \dfrac{3^{2k+2}-1}{4}\left(\dfrac{f_2^6f_3^{12}}{qf_1^6f_6^6}+20\dfrac{f_2^3f_3^3f_6^3}{f_1^3}+q\dfrac{f_6^{12}}{f_3^6}\right)\nonumber\\
	&+\dfrac{27(3^{2k}-1)(3^{2+2}-1)}{320}\left[3\left(\dfrac{f_2^6f_3^{12}}{qf_1^6f_6^6}+16\dfrac{f_2^3f_3^3f_6^3}{f_1^3}\right)+12\dfrac{f_2^3f_3^3f_6^3}{f_1^3}\right]+3^{2k+2}\cdot\dfrac{f_4^8}{f_1^2}\nonumber\\
	&= \left(\dfrac{3^{2k+2}-1}{4}+\dfrac{81(3^{2k}-1)(3^{2k+2}-1)}{320}\right)\dfrac{f_2^6f_3^{12}}{qf_1^6f_6^6}\nonumber\\
	&+\left(5(3^{2k+2}-1)+\dfrac{81(3^{2k}-1)(3^{2k+2}-1)}{16}\right)\dfrac{f_2^3f_3^3f_6^3}{f_1^3}+\dfrac{3^{2k+2}-1}{4}\cdot \dfrac{qf_6^{12}}{f_3^6}+3^{2k+2}\cdot\dfrac{f_4^8}{f_1^2}\nonumber\\
	&= \dfrac{(3^{2k+2}-1)(3^{2k+4}-1)}{320}\cdot\dfrac{f_2^6f_3^{12}}{qf_1^6f_6^6}+\dfrac{(3^{2k+2}-1)(3^{2k+4}-1)}{16}\cdot\dfrac{f_2^3f_3^3f_6^3}{f_1^3}\nonumber\\
	&+\dfrac{3^{2k+2}-1}{4}\cdot \dfrac{qf_6^{12}}{f_3^6}+3^{2k+2}\cdot\dfrac{f_4^8}{f_1^2}\nonumber\\
	&= \dfrac{(3^{2k+2}-1)(3^{2k+4}-1)}{320}\left(\dfrac{f_2^6f_3^{12}}{qf_1^6f_6^6}+20\dfrac{f_2^3f_3^3f_6^3}{f_1^3}\right)+\dfrac{3^{2k+2}-1}{4}\cdot\dfrac{qf_6^{12}}{f_3^6}+3^{2k+2}\cdot\dfrac{f_4^8}{f_1^2}.
\end{align}

We now apply $V$ on both sides of (\ref{eq42}). In view of Theorem \ref{thm31} and Lemma \ref{lem34}, we find that
\begin{align}\label{eq43}
	\sum_{n=0}^\infty &A_4\left(3^{2k+3}n+\dfrac{3^{2k+4}-5}{4}\right)q^n \nonumber\\
	&=\dfrac{(3^{2k+2}-1)(3^{2k+4}-1)}{320}\left[V\left(\dfrac{f_2^6f_3^{12}}{qf_1^6f_6^6}\right)+20V\left(\dfrac{f_2^3f_3^3f_6^3}{f_1^3}\right)\right]+\dfrac{3^{2k+2}-1}{4}\cdot\dfrac{f_2^{12}}{f_1^6}\nonumber\\
	&+3^{2k+2}\sum_{n=0}^\infty A_4(3n+1)q^n\nonumber\\
	&=\dfrac{(3^{2k+2}-1)(3^{2k+4}-1)}{320}\left[21\dfrac{f_2^2f_3^{10}}{f_1^4f_6^2}+48q\dfrac{f_3f_6^7}{f_1f_2}+20\left(3\dfrac{f_2^2f_3^{10}}{f_1^4f_6^2}+3q\dfrac{f_3f_6^7}{f_1f_2}\right)\right]\nonumber\\
	&+\dfrac{3^{2k+2}-1}{4}\cdot\dfrac{f_2^{12}}{f_1^6}+3^{2k+2}\left(2\dfrac{f_2^{12}}{f_1^6}+9q^3\dfrac{f_{12}^8}{f_3^2}\right)\nonumber\\
	&= \dfrac{3^{2k+4}-1}{4}\cdot\dfrac{f_2^{12}}{f_1^6}+\dfrac{27(3^{2k+2}-1)(3^{2k+4}-1)}{320}\left(3\dfrac{f_2^2f_3^{10}}{f_1^4f_6^2}+4q\dfrac{f_3f_6^7}{f_1f_2}\right)+3^{2k+4}\cdot q^3\dfrac{f_{12}^8}{f_3^2}.
\end{align}
We infer from (\ref{eq43}) that (\ref{eq11}) also holds for $k+1$, so by induction, (\ref{eq11}) holds for all $k\geq 0$.
\end{proof}

\begin{proof}[Proof of Theorem \ref{thm12}]
For real-valued constants $R_k$ and $S_k$ with $k\geq 0$, we examine the expansion of 
\begin{align*}
\sum_{n=0}^\infty \left[A_4\left(3^{2k+1}n+\dfrac{3^{2k+2}-5}{4}\right)-R_kA_4(27n+19)-S_kA_4(3n+1)\right]q^n.
\end{align*}
Invoking Theorems \ref{thm11}, \ref{thm31}, and \ref{thm35}, we know that the above expression equals
\begin{align}
&\dfrac{3^{2k+2}-1}{4}\cdot\dfrac{f_2^{12}}{f_1^6}+\dfrac{27(3^{2k}-1)(3^{2k+2}-1)}{320}\left(3\dfrac{f_2^2f_3^{10}}{f_1^4f_6^2}+4q\dfrac{f_3f_6^7}{f_1f_2}\right)+3^{2k+2}\cdot q^3\dfrac{f_{12}^8}{f_3^2}\nonumber\\
&-R_k\left[20\dfrac{f_2^{12}}{f_1^6}+54\left(3\dfrac{f_2^2f_3^{10}}{f_1^4f_6^2}+4q\dfrac{f_3f_6^7}{f_1f_2}\right)+81q^3\dfrac{f_{12}^8}{f_3^2}\right]\nonumber\\
&-S_k\left(2\dfrac{f_2^{12}}{f_1^6}+9q^3\dfrac{f_{12}^8}{f_3^2}\right)\nonumber\\
&= \left(\dfrac{3^{2k+2}-1}{4}-20R_k-2S_k\right)\dfrac{f_2^{12}}{f_1^6}+\left(\dfrac{27(3^{2k}-1)(3^{2k+2}-1)}{320}-54R_k\right)\left(3\dfrac{f_2^2f_3^{10}}{f_1^4f_6^2}+4q\dfrac{f_3f_6^7}{f_1f_2}\right)\nonumber\\
&+\left(3^{2k+2}-81R_k-9S_k\right)\cdot q^3\dfrac{f_{12}^8}{f_3^2}.\label{eq44}
\end{align} 
We now choose $R_k$ and $S_k$ so that the coefficients of 
\begin{align*}
\dfrac{f_2^{12}}{f_1^6}\qquad\text{ and }\qquad 3\dfrac{f_2^2f_3^{10}}{f_1^4f_6^2}+4q\dfrac{f_3f_6^7}{f_1f_2}
\end{align*}
are both zero for each $k\geq 0$. Thus, setting
\begin{align*}
\dfrac{3^{2k+2}-1}{4}-20R_k-2S_k &= 0,\\
\dfrac{27(3^{2k}-1)(3^{2k+2}-1)}{320}-54R_k &=0,
\end{align*}
and solving this system of linear equations, we get
\begin{align*}
(R_k,S_k) = \left(\dfrac{(3^{2k}-1)(3^{2k+2}-1)}{640},-\dfrac{(3^{2k}-9)(3^{2k+2}-1)}{64}\right).
\end{align*}
Observe that for all integers $k\geq 0$, $3^{2k}-1, 3^{2k}-9$, and $3^{2k+2}-1$ are all divisible by $8$, and exactly one of $3^{2k}-1$ and $3^{2k+2}-1$ is divisible by $10$. Thus, we see that $R_k$ and $S_k$ are integers for all integers $k\geq 0$. We also deduce from (\ref{eq44}) that
\begin{align}
\sum_{n=0}^\infty &A_4\left(3^{2k+1}n+\dfrac{3^{2k+2}-5}{4}\right)q^n\nonumber\\
&=\dfrac{(3^{2k}-1)(3^{2k+2}-1)}{640}\sum_{n=0}^\infty A_4(27n+19)q^n-\dfrac{(3^{2k}-9)(3^{2k+2}-1)}{64}\sum_{n=0}^\infty A_4(3n+1)q^n\nonumber\\
&+\dfrac{81(3^{2k}-1)(3^{2k}-9)}{640}\cdot q^3\dfrac{f_{12}^8}{f_3^2}.\label{eq45}
\end{align}
Looking at the terms involving $q^{3n+1}$ and $q^{3n+2}$ in the expansion of (\ref{eq45}), we find that 
\begin{align}\label{eq46}
	A_4\left(3^{2k+1}n+\dfrac{3^{2k+2}-5}{4}\right) &= \dfrac{(3^{2k}-1)(3^{2k+2}-1)}{640}A_4(27n+19)\nonumber\\
	&\phantom{test}-\dfrac{(3^{2k}-9)(3^{2k+2}-1)}{64}A_4(3n+1)
\end{align}
for all integers $n\geq 0$ with $3\nmid n$ and $k\geq 0$. Combining (\ref{eq46}) and Theorem \ref{thm36} yields (\ref{eq12}). Hence, (\ref{eq13}) now follows from (\ref{eq12}) and Theorem \ref{thm32}.
\end{proof}

We close this section by presenting another linear identity for $A_4(n)$ and give some remarks.

\begin{theorem}\label{thm41}
For all integers $n\geq 0$ and $k\geq 0$, we have
\begin{align}\label{eq47}
A_4\left(3^{2k+2}n+\dfrac{5(3^{2k+2}-1)}{4}\right) &= \dfrac{(3^{2k}-1)(3^{2k+2}-1)}{640}A_4(81n+100)\nonumber\\
&-\dfrac{(3^{2k}-9)(3^{2k+2}-1)}{64}A_4(9n+10)+\dfrac{81(3^{2k}-1)(3^{2k}-9)}{640}A_4(n). 
\end{align}
\end{theorem}

\begin{proof}
Dividing both sides of (\ref{eq45}) by $q^3$, we have 
\begin{align}
	\sum_{n=-3}^\infty &A_4\left(3^{2k+1}n+\dfrac{5(3^{2k+2}-1)}{4}\right)q^n\nonumber\\
	&=\dfrac{(3^{2k}-1)(3^{2k+2}-1)}{640}\sum_{n=-3}^\infty A_4(27n+100)q^n-\dfrac{(3^{2k}-9)(3^{2k+2}-1)}{64}\sum_{n=-3}^\infty A_4(3n+10)q^n\nonumber\\
	&+\dfrac{81(3^{2k}-1)(3^{2k}-9)}{640}\cdot\dfrac{f_{12}^8}{f_3^2}.\label{eq48}
\end{align}
We now apply $U$ on both sides of (\ref{eq48}) and use the generating function for $A_4(n)$, so that
\begin{align}
	\sum_{n=-1}^\infty &A_4\left(3^{2k+2}n+\dfrac{5(3^{2k+2}-1)}{4}\right)q^n\nonumber\\
	&=\dfrac{(3^{2k}-1)(3^{2k+2}-1)}{640}\sum_{n=-1}^\infty A_4(81n+100)q^n-\dfrac{(3^{2k}-9)(3^{2k+2}-1)}{64}\sum_{n=-1}^\infty A_4(9n+10)q^n\nonumber\\
	&+\dfrac{81(3^{2k}-1)(3^{2k}-9)}{640}\sum_{n=0}^\infty A_4(n)q^n.\label{eq49}
\end{align}
Hence, (\ref{eq47}) follows from comparing the coefficients of $q^n$ with $n\geq 0$ in the expansion of (\ref{eq49}).
\end{proof}

\begin{remark}\label{rem42}
As mentioned in the proof of Theorem \ref{thm12}, we know that for all integers $k\geq 0$, $3^{2k}-1, 3^{2k}-9$, and $3^{2k+2}-1$ are all multiples of $8$. Also, exactly one of $3^{2k}-1$ and $3^{2k+2}-1$ is divisible by $10$, and exactly one of $3^{2k}-1$ and $3^{2k}-9$ is divisible by $10$. Since 
\begin{align*}
\dfrac{(3^{2k}-1)(3^{2k+2}-1)}{640}-\dfrac{(3^{2k}-9)(3^{2k+2}-1)}{64}+\dfrac{81(3^{2k}-1)(3^{2k}-9)}{640}=1
\end{align*}
for all integers $k\geq 0$, we surmise that
\begin{align*}
\dfrac{(3^{2k}-1)(3^{2k+2}-1)}{640}, \qquad \dfrac{(3^{2k}-9)(3^{2k+2}-1)}{64}, \qquad \dfrac{81(3^{2k}-1)(3^{2k}-9)}{640}
\end{align*}
are all integers having no common factor. This shows that there are no nontrivial (infinite family of) congruences for $A_4(n)$ such as (\ref{eq13}) that can be found via (\ref{eq47}). 
\end{remark}

\bibliography{bipfourcores}
\end{document}